\documentclass[12pt]{article} 
\usepackage[dvips]{graphics,color}

\usepackage{latexsym}
\usepackage{amssymb}

\usepackage{latexsym}
\usepackage{amssymb}
\usepackage{euscript}
\usepackage[dvips]{graphics}
\usepackage{epsf}

\def\mcc{M\raise.5ex\hbox{c}C}
\def\mccarthy{M\raise.5ex\hbox{c}Carthy}

\def\ie{{\it i.e. }}

\def\a{\alpha}

\def\={\ = \ }

\def\C{\mathbb C}
\def\R{\mathbb R}

\def\D{\mathbb D}

\def\inn{\ \in \ }

\def\be{\setcounter{equation}{\value{theorem}} \begin{equation}}
\def\ee{\end{equation} \addtocounter{theorem}{1}}
\def\beq{\begin{eqnarray*}}
\def\eeq{\end{eqnarray*}}

\def\bp{{\sc Proof: }}
\def\ep{{}{\hfill $\Box$} \vskip 5pt \par}

\def\bl{\begin{lemma}}
\def\el{\end{lemma}}
\def\bt{\begin{theorem}}
\def\et{\end{theorem}}
\def\bprop{\begin{prop}}
\def\eprop{\end{prop}}
\def\bd{\begin{definition}}
\def\ed{\end{definition}}
\def\br{\begin{remark}}
\def\er{\end{remark}}
\def\bexer{\begin{exercise}}
\def\eexer{\end{exercise}}
\def\bfig{\begin{figure}}
\def\efig{\end{figure}}

\newtheorem{theorem}{Theorem}[section]
\newtheorem{prop}[theorem]{Proposition}
\newtheorem{lemma}[theorem]{Lemma}
\newtheorem{cor}[theorem]{Corollary}

\newtheorem{definition}[theorem]{Definition}

\newcommand{\al}{\alpha}

\begin{document}
\setlength{\baselineskip}{21pt}
\title{Cusp algebras}
\author{Jim Agler
\thanks{Partially supported by National Science Foundation Grant
DMS 0400826}\\
U.C. San Diego\\
La Jolla, California 92093
\and
John E. M\raise.5ex\hbox{c}Carthy
\thanks{Partially supported by National Science Foundation Grant
DMS 0501079}\\
Washington University\\
St. Louis, Missouri 63130}

\bibliographystyle{plain}

\maketitle
\begin{abstract}
\end{abstract}

\baselineskip = 18pt

\section{Introduction}\label{seca}

By a {\em cusp} $V$ we shall mean the image of the unit disk $\D$ under a bounded injective holomorphic 
map $h$ into $\C^n$ whose derivative vanishes at exactly one point.
The simplest example is the Neil parabola, given by $h (\zeta)
= (\zeta^2,\zeta^3)$.
See \cite{kn07a, pausin06}
for background and theory on the Neil parabola, which is pictured in Figure~\ref{figa}.

A generalization of a cusp is a {\em petal}.
A petal is
the image of the unit disk $\D$  under a proper holomorphic map 
$h$ from $\D$ into some bounded open set $\Omega$,
where, 
except for a finite set $E_h$,
$h$ is one-to-one and non-singular. 
If $E_h$ is a singleton, then $V$ is  a cusp.
Since the automorphism group of $\D$ is transitive, we may assume then that $E_h = \{ 0 \}$. 
\bfig 
\includegraphics{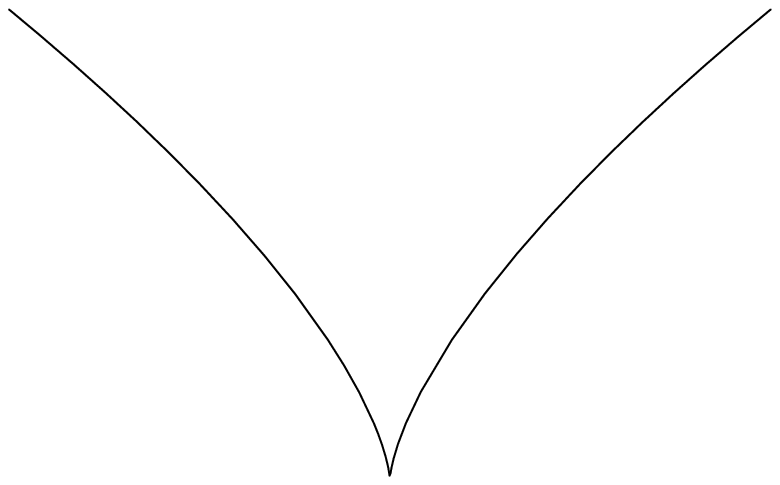} \caption{A simple cusp - the Neil parabola $z^3 = w^2$} 
\label{figa} 
\efig

The function $h :\ \D \to \Omega $ induces a finite codimensional subalgebra $A_h$ of the algebra
$O(\D)$ of all analytic functions on $\D$, namely 
$$
A_h \= \{ F \circ h \ : \ F \inn O(V) \},
$$
where $O(V)$ is the algebra of functions on $V$ that, in the neighborhood of every point
$P$ of $V$, coincides with the restriction to $V$ of a 
function holomorphic in a neighborhood of $\C^n$ containing $P$.
When $V$ is a cusp, the algebra $A_h$ has the property that it contains
$z^m O(\D)$ for some $m \geq 2$. Conversely, if a cofinite subalgebra 
of $O(\D)$ contains $z^m O(\D)$, then it arises in this way.

In light of these remarks, let us agree to say that a
unital  algebra $A \subseteq O(\D)$ is a cusp algebra
if

(i) dim$\left( O (\D) / A \right)$ is finite, and

(ii) For some $m \geq 2$, $z^m O(\D) \subseteq A$.

We shall say a map $h: \D \to \Omega \subsetneq \C^n$ is a {\em holomap} if it
is
a proper map that is one-to-one and non-singular except on a finite set.
We shall say that a holomap $h$ is a 
{\em holization} of the cofinite algebra $A$ if $A = A_h$.

If $A$ is a cusp algebra, then we define three integers attached to $A$.
The {\em codimension} of $A$ is defined by
$$
{\rm cod}(A) \= {\rm dim}\left( O(\D) / A \right) .
$$
The {\em order} of $A$
is defined by
$$
{\rm ord}(A) \= \min\{ k \ : \ z^{k+1} O(\D) \subseteq A \} .
$$
The {\em contact} of $A$ is
defined by
$$
{\rm con}(A) \= \max\{ n \ : \ f^{(j)} (0) = 0, \ \forall \ f \inn A, \ 1 \leq j \leq n \} .
$$
If $V$ is a cusp, then we define ${\rm cod}(V)$, ${\rm ord}(V)$
and
${\rm con}(V)$
as the quantities for the corresponding algebra $A_h$. 

By a {\em simple} cusp algebra, we shall mean a cusp algebra $A$
such that ${\rm con}(A) = 1$. There is a unique cusp algebra of codimension
$1$, which corresponds via a holization to the holomorphic functions on the 
Neil parabola $\{ (z,w) \inn \D^2: z^3 = w^2 \}$ (this result can be seen directly, or
will follow from Section~\ref{secd}.) 
This can be thought of as a Riemann mapping theorem for cusps of codimension $1$.
In Section~\ref{secd}, we generalize this ``Riemann mapping theorem'' to arbitrary simple cusps.
We show that locally the only invariant of simple cusps is the codimension, but globally
cusps of codimension $n+1$ have a $2n-1$ real parameter moduli space.

This paper is a continuation of the authors' earlier work in \cite{amhyp}, where the ideas 
are put in a more general context.
A principal concern in that paper is when a given finite codimensional subalgebra 
can be holized in 2 dimensions, \ie by a map into $\C^2$. In Section~\ref{sece}, we show that all
simple cusp algebras can be holized in $\C^2$.

\section{Preliminaries}\label{secb}

Throughout this section, we shall assume that $A$ is a simple cusp algebra, \ie of contact one.
By a {\em primitive} for $A$ we mean a function $\pi$ that lies in $A$, and satisfies
$\pi(0) = 0,\ \pi''(0) = 2$; so $\pi $ has a Taylor expansion
$\pi(z) = z^2 + \dots$.

For each $n$, let $A_n$ be the space
$$
A_n \= \left[ A \cap z^{2n} O(\D) \right],
$$ and let $E_n$ be the quotient
$$
E_n \= A_n / A_{n+1} .
$$
Observe that multiplication by any primitive $\pi$ is always one-to-one from $E_n$ to $E_{n+1}$.
Therefore, there is some integer $n_0$ such that 
${\rm dim}(E_n)$ is $1$ for $n \leq n_0$, and $2$ for $n > n_0$.
Since ${\rm dim}(E_n) = 2$
if and only if both $z^{2n}$ and $z^{2n+1}$ are in $A$, it follows that 
\beq
\label{eqb1}
{\rm ord}(A) \= 2n_0 + 1 .
\eeq
As $\pi^n$ is in $A_n$ for each $n$, if $f$ is in $A_n$ for some $n \leq n_0$,
then there is a constant $c_n$ such that $f - c_n \pi^n$ is in $A_{n+1}$.
Therefore every function $f$ in $A$ has the representation
\be
\label{eqb2}
f \= c_0 + c_1 \pi + \dots + c_{n_0} \pi^{n_0} + z^{2n_0 +2} g ,
\ee
where $c_0,\dots,c_{n_0}$ are in $\C$ and $g \inn O(\D)$.
From (\ref{eqb2}), we see that ${\rm cod}(A) = n_0 + 1$.
Summarizing, we have:
\bprop
\label{propb1}
If $A$ is a simple cusp algebra, then ${\rm ord}(A) = 2\,{\rm cod}(A) -1$.
Furthermore, if we set $n_0 = {\rm cod}(A) -1$, and $\pi$ is a primitive for $A$,
then every function $f$ in $A$ has a unique representation
\be
\label{eqb3}
f(z) \= p \circ \pi (z) \ + \ z^{2n_0 +2} g(z) 
\ee
for some polynomial $p$ of degree at most $n_0$ and some $g \inn O(\D)$.
\eprop

\section{Connections and local theory}\label{secc}

Let $U$ be an open set in $\C$.
A linear functional on $O(U)$ is called local if it comes from a
finitely supported distribution, \ie is of the form
$$
\Lambda(f) \= \sum_{i=1}^m \sum_{j=0}^{n_i} a_{ij}
f^{(j)}(\alpha_i) .
$$

\bd
\label{def22}
A connection on $\{ \a_1,\dots , \a_m \} \subset U$
is a finite dimensional set $\Gamma$ of local functionals $\Lambda$ supported by
$\{ \a_1,\dots , \a_m \}$. We say $\Gamma$ is algebraic if
$\Gamma^\perp := \{ f \ \in \ O(U)\ : \ \Lambda(f) = 0 \ \forall \,
\Lambda\, \in \, \Gamma \}$ is an algebra.
\ed

It was proved by T. Gamelin \cite{gam68} that every finite codimensional subalgebra $A$ of $O(\D)$ is
$\Gamma^\perp$ for some algebraic connection on $\D$. Moreover, $A$ will be a cusp algebra 
iff the support of $\Gamma$ is $\{ 0 \}$.

We shall say that  a point $P$ in a petal $V$ is a {\em cusp point} if there is a one-to-one proper map $h$ 
from $\D$ onto a neigborhood $U$ of $P$ in $V$ such that the derivative of $h$ vanishes only
at the pre-image of $P$. The algebra $A_h$ is then a cusp algebra.
We wish to show that its codimension, order and contact do not depend on the choice of $U$.

Suppose that $P_1$ and $P_2$ are cusp points in cusps $V_1$ and $V_2$, and there is an
injective holomorphic map $\phi: V_1 \to V_2$. 
For $r=1,2$ there are holizations $h_r : \D \to V_r$ with $h_r(0) =P_r$.
By \cite[Thm. 3.2]{amhyp}, the map $\phi$ from $V_1$ to $V_2$ induces a 
one-to-one map $\psi : \D \to \D$
such that $ h_2 \circ \psi = \phi \circ h_1$,
and the connections $\Gamma_1$ and $\Gamma_2$ induced
by $h_1$ and $h_2$ are related by 
\beq
\Gamma_2 &\ \supseteq \ & \psi_{*}\, \Gamma_1 \\
&\ :=\ & \{ f \mapsto \Lambda_1 (f \circ \psi) \ : \ \Lambda_1 \inn \Gamma_1 \}.
\eeq

If 
\be
\label{eqc01}
\Lambda_1 \ : \ f \mapsto \sum_{j=0}^n a_j  f^{(j)}(0),
\ee
then
\be
\label{eqc1}
\psi_*\, \Lambda_1 \ : \ g \mapsto \sum_{j=0}^n a_j  (g\circ \psi)^{(j)}(0),
\ee
We can use Fa\'a di Bruno's formula to evaluate the derivatives of
$g\circ\psi$ in terms of those of $g$ and $\psi$.
As $\psi(0) = 0$ and $\psi'(0) \neq 0$,
note in particular that the coefficient of $g^{(n)}(0)$ on the right-hand side of
(\ref{eqc1}) is $a_n [\psi'(0)]^n$. 
Therefore $\psi_*$ is injective, so ${\rm cod}(A_1) \leq {\rm cod}(A_2)$.

If ${\rm con}(A_1) = k$, it means that for each $1 \leq j \leq k$,
the functional $f \mapsto f^{(j)}(0)$ is in $\Gamma_1$. Therefore 
each 
functional $g \mapsto g^{(j)}(0)$ is in $\Gamma_2$, and
${\rm con}(A_1) \leq {\rm con}(A_2)$.
Finally if 
${\rm ord}(A_1) = n$, it means there is 
some $\Lambda_1$ in $\Gamma_1$ of the form (\ref{eqc01}) with $a_n \neq 0$;
it follows that ${\rm ord}( A_1) \leq {\rm ord}(A_2)$.

Summarizing, we have shown:
\bprop
\label{propc1}
If $\phi : V_1 \to V_2$ is an injective holomorphic map from one cusp to another, then
${\rm cod}(V_1) \leq {\rm cod}(V_2),\ {\rm con}(V_1) \leq {\rm con}(V_2)$ and
${\rm ord}(V_1) \leq {\rm ord}(V_2)$.
\eprop

Now, if $P$ is a cusp point in a petal $V_1$, and one surrounds $P$ by a decreasing sequence
of open sets $U_n$, then with each $U_n$ there is a cusp algebra with a codimension, contact and order,
and there is an inclusion map from each $U_{n+1}$ to $U_n$.
By Proposition~\ref{propc1}, the positive integers 
${\rm cod}(U_n), \ {\rm con}(U_n)$ and
${\rm ord}(U_n)$
must decrease as $U_n$ shrinks. Therefore, 
at some point they must stabilize, and for all sufficiently small neighborhoods of $P$, the cusp algebras
have the same codimension, contact and order. We shall call these ${\rm }(cod)(P),\ {\rm }(con)(P)$
and ${\rm }(ord)(P) $ respectively, and say $P$ is a simple cusp point if ${\rm con}(P) = 1$. 

\section{Equivalence of simple cusps}\label{secd}

We shall say that two points $P_1$ and $P_2$ in the
petals $V_1$ and $V_2$ are 
{\em locally equivalent} if there is a  neighborhood $U_1$ of $P_1$ in $V_1$,
and a neighborhood $U_2$ of $P_2$ in $V_2$,
and a biholomorphic homeomorphism $F$ from $U_1$ onto $U_2$ that maps $P_1$ to $P_2$.

\bt
\label{thmd1}
Two simple cusp points are locally equivalent if and only 
they have the same codimension.
\et
\bp
Necessity follows from Proposition~\ref{propc1}. 
To prove sufficiency, 
let $\pi_1$ be a primitive of the first algebra, and $\pi_2$ a primitive
of the second algebra.
By Proposition~\ref{propb1},
it is sufficient to prove that there are neighborhoods $W_1$ and $W_2$ of the origin,
and a univalent map $\phi : W_1 \to W_2$  that maps $0$ to $0$ and
such that
\be
\label{eqd8}
\pi_2 \circ \phi \= \pi_1 .
\ee
Each $\pi_r$ has a square root $\chi_r$ in a neighborhood of $0$, and each $\chi_r$
is locally univalent (because $\pi_r$ is of order $2$). 
So define 
$$
\phi \ := \ \chi_2^{-1} \circ \chi_1
$$
on a suitable neighborhood $W_1$, and (\ref{eqd8}) is satisfied.
\ep

How can we globally parametrize simple cusp algebras?
We shall show that there is an essentially unique primitive with all
its even Taylor coefficients (except for the second) zero. By
$\hat \pi (k)$ we mean the $k^{\rm th}$ Taylor coefficient at $0$.
\bl
\label{lemd1}
Every simple cusp algebra of codimension $n+1$ has a primitive $\pi$
such that $\hat \pi (2k) \= 0$ for all $ 2 \leq k \leq n$.
Moreover, $\pi$ is unique up to $O(z^{2n+2})$.
\el
\bp
(Existence)
Let $\chi_1$ be any primitive. Define
$$
\chi_2 \ := \ \chi_1 - \hat \chi_1(4) (\chi_1)^2 .
$$
Proceed inductively, with
$$
\chi_k \ := \ \chi_{k-1} - \hat \chi_{k-1}(2k) (\chi_{k-1})^k .
$$
Then let $\pi = \chi_n$.

(Uniqueness)
Suppose $\pi$ and $\chi$
are both primitives with their even Taylor coefficients, starting at 4, vanishing.
By Proposition~\ref{propb1}, we have
$$
\chi \= c_1 \pi + c_2 \pi^2 + \dots + c_n \pi^n + O(z^{2n+2}) .
$$
By looking at the coefficient of $z^2$ we see $c_1 = 1$.
Now looking at the coefficients of $z^4, z^6,\dots,z^{2n}$ in order, we see that
$c_2 = 0 = c_4 = \dots = c_n$.
\ep
So every algebraic connection of dimension $n+1$ supported at the origin is the annihilator 
of an algebra generated locally by a unique primitive
\be
\label{eqd12}
\pi(z) \= z^2 + \alpha_1 z^3 + \alpha_2 z^5 + \dots + \alpha_n z^{2n+1} .
\ee
Let us denote the 
codimension $n+1$ subalgebra of $O(\D)$ with primitive (\ref{eqd12})
by $A(\a_1,\dots,\a_n)$,
and let $V(\a_1,\dots,\a_n)$ denote the corresponding petal 
$h(\D)$, where $h$ holizes $A(\a_1,\dots,\a_n)$.

Two cusps $V_1$ and $V_2$ are globally equivalent if there is a biholomorphic homeomorphism $\phi$
from
$V_1$ onto $V_2$. If each $V_r$ is holized by $h_r : \D \to V_r$,
then by \cite[Thm. 3.1]{amhyp}, this occurs if and only if there is a map $\psi: \D \to \D$ 
such that
$\psi_*$ maps the first connection to the second. As the only self-maps of the disk that leave the origin 
invariant are rotations, this is rather restrictive. 

\bt
\label{thmd2}
The cusps $V(\a_1,\dots,\a_n)$ and $V(\beta_1,\dots,\beta_n)$
are isomporphic if and only if there is a unimodular constant $\tau$ such that
$\beta_j = \tau^{2j-1} \a_j$ for $1 \leq j \leq n$.
\et
As a corollary we have
\begin{cor}
\label{cord1}
The moduli space of all simple cusps of codimension $n+1$ is $\R^+ \times \C^{n-1}$.
\end{cor}

\section{Embedding}\label{sece}

The purpose of this section 
is to prove that every simple cusp algebra $A$ contains a pair of functions
$h_1,h_2$ such that the pair holizes the algebra (which is equivalent to saying
that polynomials in $h_1$ and $h_2$
are dense in $A$ in the topology of uniform convergence on compacta).

For the rest of  this section, $A$ will be a fixed simple cusp algebra of codimension $n+1$.
\bl
\label{leme1}
For every $\alpha$ on $\D \setminus \{0\}$, there is a function $\psi_\alpha$ in
$A$ that has a single simple zero at $\alpha$, and no other zeroes on $\D$.
\el
\bp
Consider functions of the form $(z-\al)e^{h(z)}$. In order to be in $A$,
$h$ must satisfy $n+1$ equations on its first $2n+1$ derivatives at $0$.
This system is triangular, so can be solved with $h$ a polynomial.
\ep
\bl
\label{leme2}
If $f$ is in $A$ and $f$ has no zeroes on $\D$, then $1/f$ is in $A$.
\el
\bp
Consider the vector space obtained by adjoining $1/f, 1/f^2,\dots, 1/f^k$ to $A$.
Since $A$ is of finite codimension in $O(\D)$, for some $k$ there is a linear relation
$$
1/f^k \= \sum_{j=1}^{k-1} c_j 1/f^j \ + \ A .
$$
Multiplying both sides by $f^{k-1}$, we get that $1/f$ is in $A$.
\ep
\bl
\label{leme3}
If $f \inn A$ and $f(\al) = 0$ for some $\alpha \inn \D \setminus \{ 0 \}$, then 
$f/\psi_\al$ is in $A$.
\el
\bp
Let $\Gamma$ be the connection at $0$ such that $A = \Gamma^\perp$.
Consider the algebra 
$$
A' \ := \ \{ f \inn O (\frac{|\al|}{2} \D) \ : \ \Lambda(f) = 0 \ \forall\  \Lambda \inn \Gamma \} .
$$
Then $\psi_\al$ is in $A'$, so by Lemma~\ref{leme2}, $1/\psi_\al$ is in $A'$.
Therefore $f/\psi_\al$ is in $A'$, and is therefore annihilated by every functional in $\Gamma$.
But $f/\psi_\al$ is also in $O(\D)$, so it is in $A$ as required.
\ep
\bt
\label{thme1}
There exists a map $h: \D \to \C^2$ that holizes $A$.
\et
\bp
By applying Lemma~\ref{leme3} repeatedly, one can find a primitive $h_1$ of $A$ that
has no zeroes in $ \D \setminus \{ 0 \}$.
Define
$$
h_2(z) \ := \ z \, (h_1(z) )^{n + 1} .
$$
Consider the algebra $A'$ that is the closure (in the topology of uniform convergence on 
compacta) of polynomials in $h_1$ and $h_2$.
By \cite[Thm. 4.2]{amhyp}, the algebra $A'$ is of finite codimension, and it 
is a cusp algebra, because $h(z) = (h_1(z),h_2(z))$ is one-to-one away from the origin.
If $A' \neq A$, then it is contained in a maximal proper subalgebra $A''$ of $A$, which must
have codimension $n+2$, and therefore order $2n+3$.
But $A''$ contains $h_2$ and $h_1^{n+1}$, and these functions have order $2n+2$ and
$2n+3$ respectively
at the origin. So there can be no linear relation in $A''$ between the derivatives at the origin
of order $(2n+2)$
and $(2n+3)$, and so the order of $A''$ must actually be at most $2n +1$.
\ep
\bibliography{references}

\end{document}